\documentclass[A4paper,12pt]{book}
\usepackage{amssymb,amsmath}
\usepackage[russian]{babel}
\usepackage{graphicx}
\newcount\exnom\exnom=0

\long\def\exo/{\vspace{0.1cm} \noindent\advance\exnom by1{
{\the\exnom}.}}

\begin{document}
\thispagestyle{empty}

\begin{center} \textbf{DERIVATIVES OF TRIANGULAR, TOEPLITZ, CIRCULANT MATRICES AND OF MATRICES OF OTHER FORMS OVER SEMIRINGS}
\end{center}

\begin{center}
\textbf{Dimitrinka Vladeva}\\
Department of  Mathematics and Physics,\\
  Forest University, 10, Kliment Ohridski Str.\\
1756 Sofia, Bulgaria
\end{center}

 \vspace{2mm}

\noindent\emph{\textbf{Abstract:} In this article we construct  examples of derivations in matrix semirings. We study hereditary and inner derivations, derivatives of diagonal, triangular, Toeplitz, circulant matrices and of matrices of other forms and prove theorems for derivatives of matrices of these  forms.}


\vspace{5mm}

\begin{center} \textbf{Introduction and preliminaries}
\end{center}

\vspace{3mm}

 An algebra $\;S = (S,+,.)\;$ with two binary operations $\;+\;$ and $\;\cdot\;$ on $\;S$, is called a \textbf{{{semiring}}} if:

\textbf{1.} $\;(S,+)$ is a commutative semigroup,

\textbf{2.} $\; (S,\cdot)$ is a semigroup,

\textbf{3.}   distributive laws hold $ x\cdot(y + z) = x\cdot y + x\cdot z$ and $(x + y)\cdot z = x\cdot z + y\cdot z$
for any $x, y, z \in S$.
\vspace{1mm}

The semiring $S$ is called \textbf{commutative} if $a\cdot b = b\cdot a$ for any $a, b \in S$.
\vspace{1mm}

In all considered semirings $S$ we assume that there exist a \textbf{zero} element $0 \in S$ such that $a + 0 = a$ and $a\cdot 0 = 0\cdot a = 0$ for any $a \in S$ and \textbf{identity} element $1 \in S$ such that $1\cdot a = a\cdot 1 = a$ for any $a \in S$.
\vspace{2mm}

An element $a$ of a semiring $S$ is called \textbf{additively  idempotent} if $a + a = a\;$.
A semiring $S$ is called \textbf{{additively idempotent}} if each of its elements is additively idempotent. Additively idempotent semirings are proper semirings, i.e. they are not rings.  Basic facts for semirings  can be found in  [6]. For importance of matrices over additively idempotent semirings and semifields and their applications the reader can see [2], [4] and [9]. Combinatorial ideas in matrix theory was obtained in [3].
\vspace{3mm}

Let $S$ be a semiring. A map $\delta : S \rightarrow S$ is a \textbf{derivation} of $S$ if and only if
$$\delta(a + b) = \delta(a) + \delta(b) \; \; \mbox{and}\; \; \delta(a\,b) = \delta(a)\,b + a\,\delta(b),$$
where $a, b \in S$.

Note that identity map is a derivation in any additively idempotent semiring.

 An overview of derivation in semirings can be found in [5].
\vspace{3mm}

In the present article we study the derivations in matrix semirings. Let $S_0$ be a semiring and  $M_n(S_0)$ be the semiring of $n\times n$ - matrices with entries from $S_0$. Let $S$ be a subsemiring of $M_n(S_0)$, containing the zero matrix and identity matrix $E$. A map $\delta : S \rightarrow S$ is a derivation if for any matrices $A, B \in S$ and $\alpha \in S_0$, it follows
$$\begin{array}{l} (1) \;\; \delta(A + B) = \delta(A) + \delta(B),\\
(2)\;\; \delta(A\,B) = \delta(A)\,B + A\,\delta(B),\\
(3)\;\; \delta(\alpha A) = \delta((\alpha E) A) = \delta(\alpha E) A + \alpha \delta(A)\\
\end{array}.$$

Clearly the equality (3) is a consequence of (2). In some of constructions of derivations in the next pages $\delta(\alpha E) = 0$, but in another (where $S_0$ is an additively idempotent semiring) $\delta(\alpha E) A = \alpha \delta(A)$.

\vspace{2mm}

In [10] there is a  construction of a map from endomorphism semiring on its arbitrary subsemiring which is a projection. This projection is a derivation $\delta$ and $\delta^2 = \delta$, so it is an \textbf{idempotent derivation}.

If we extend the set of natural numbers $\mathbb{N}$ with $0$, the set $\mathbb{N}_0 = \mathbb{N}\cup\{0\}$ is a commutative semiring with zero and identity element. In the  semirig $\mathbb{N}_0[x]$, consisting of polynomials with coefficients in $\mathbb{N}_0$, one can define a derivation $\delta$ by well-known rule
$$\delta(a_0 + a_1x + \cdots + a_kx^k) = a_1 + \cdots + ka_k\,x^{k-1}.$$
Since $\delta^{k+1}(a_0 + a_1x + \cdots + a_kx^k) = 0$, it follows that  $\delta^{k+1} = 0$, so $\delta$ is a \textbf{nilpotent derivation}.
\vspace{1mm}

Some of examples of derivations in matrix semirings considered here are idempotent or nilpotent derivations.

\vspace{5mm}

\begin{center} \textbf{Hereditary derivations}
\end{center}

\vspace{3mm}

\textbf{Proposition 1.} \textsl{Let $S_0$ be a semiring and $\delta_0 : S_0 \rightarrow S_0$ be a derivation. Let $M_n(S_0)$ be a semiring of $n\times n$ - matrices  with entries from $S_0$. For $A = (a_{ij}) \in M_n(S_0)$ let $\delta_0^h(A) = \left(\delta_0(a_{ij})\right)$. Then $\delta_0^h : M_n(S) \rightarrow M_n(S)$ is a derivation.}
\vspace{1mm}

\emph{Proof.} Let $B = (b_{ij}) \in M_n(S_0)$. Then $\delta_0^h(B) = \left(\delta_0(b_{ij})\right)$ and $\delta_0^h(A + B) = \left(\delta_0(a_{ij} + b_{ij})\right) = \left(\delta_0(a_{ij}) + \delta_0(b_{ij})\right) = \left(\delta_0(a_{ij})\right) + \left(\delta_0(b_{ij})\right) = \delta_0^h(A) + \delta_0^h(B)$. Now we compute
 $$\delta_0^h(A\,B) = \left(\delta_0(\sum_{k=1}^n a_{ik} b_{kj})\right) =  \left(\sum_{k=1}^n \delta_0(a_{ik} b_{kj})\right) = \left(\sum_{k=1}^n \delta_0(a_{ik}) b_{kj} + a_{ik}\delta_0(b_{kj})\right) =$$ $$= \left(\sum_{k=1}^n \delta_0(a_{ik}) b_{kj}\right) + \left(\sum_{k=1}^n a_{ik}\delta_0(b_{kj})\right) = \delta_0^h(A)\,B + A\,\delta_0^h(B).$$

 Here the equality (3) from the previous section can be expressed directly. For $\alpha \in S_0$, it follows
 $$\delta_0^h(\alpha A) = \delta_0^h(\alpha (a_{ij})) = \left(\delta_0(\alpha a_{ij})\right) = \left(\delta_0(\alpha) a_{ij}\right) + \left(\alpha \delta_0( a_{ij})\right) = \delta_0(\alpha) A + \alpha \delta_0^h(A).$$

Derivations of this type are called \textbf{hereditary} derivations. It is easy to see that hereditary derivation $\delta_0^h$ is idempotent or nilpotent derivation if and only if the derivation $\delta_0$ is idempotent or nilpotent.
\vspace{1mm}

\textbf{\emph{Remark.}}  For hereditary derivation $\delta_0^h$, it follows
$$\delta_0^h(E)= \left(
    \begin{array}{cccc} \vspace{1mm}
      \delta_0(1) & 0 & \cdots & 0 \\ \vspace{1mm}
      \!0 & \delta_0(1) & \cdots & 0 \\
     \!\! \vdots &\!\! \vdots &\!\! \ddots & \vdots \\
     \! 0 & \!0 & \cdots & \delta_0(1) \\
    \end{array}
  \right) = \delta_0(1) E.$$

In general case the last matrix is not the zero matrix.
For example, see [10], let $S$ be a endomorphism semiring and $\delta_0 : S \rightarrow S$ be a projection on the least string of $S$, i.e. the subsemiring consisting of endomorphisms with images in set $\{0,1\}$. Since unit element of $S$ is an identity map $i$, then $\delta_0(i) = \alpha$ and $\alpha \neq 0$ since $\alpha(0) = 0$ and $\alpha(k) = 1$ for arbitrary $k > 0$.
\vspace{1mm}

Important applications of hereditary derivations can be found in [8].

\vspace{5mm}

\begin{center}
 \textbf{ Two examples of idempotent derivations}
\end{center}

\vspace{3mm}

Here we explore two semirings of matrices not considered in the next sections.
\vspace{1mm}

\textbf{\emph{Example 1. }} Consider the $n \times n$ - matrices  with entries from arbitrary (not necessary additively idempotent) semiring of the type
$$A = \left(
    \begin{array}{ccccc}\vspace{1mm}
      a_{11} & 0 & \cdots & 0 & a_{1n} \\\vspace{3mm}
     \! 0 & a_{22} & \cdots & 0 & 0\\ \vspace{3mm}
      \!\!\vdots & \!\!\!\vdots & \!\!\!\vdots & \!\!\!\!\ddots & \vdots\\
     \!\! 0 &\!\! 0 &\!\! \cdots &\! 0 & \,a_{nn} \\
    \end{array}
  \right).$$ Clearly if $B = (b_{ij})$ is a matrix of the same type then $A + B$ and $AB$ are also of the same type.

 In the semiring of matrices of this type we define a map $\delta$ such that $$\delta(A) = a_{1n}E_{1n}.$$

 Here $E_{ij}$ are the matrix units, i.e. $E_{ij} = (a_{pq})$, where $a_{ij} = 1$ and $a_{pq} = 0$ for $p \neq i$ or $q \neq j$.

Then $\delta(B) = b_{1n}E_{1n}$, $\delta(A\,B) = (a_{11}b_{1n} + a_{1n}b_{nn})E_{1n}$.

Since $\delta(A)\,B = a_{1n}b_{nn}E_{1n}$ and $A\,\delta(B) = a_{11}b_{1n}E_{1n}$, it follows that
$\delta(A\,B) = \delta(A)\,B + A\,\delta(B)$. Obviously $\delta(A + B) = \delta(A) + \delta(B)$. Hence the map $\delta$ is a derivation.
\vspace{1mm}

Since $\delta^2(A) = \delta(A)$ for any matrix $A$ of the type considered, then $\delta$ is an idempotent derivation.

 \vspace{5mm}

\textbf{\emph{Example 2. }}  Let positive integers $i_1, \ldots, i_k$ are fixed, where $k < n$ and  $S_0$ be an arbitrary semiring.

Let us consider the matrices $A = (a_{ij})$, $i, j = 1, \ldots n$, where $a_{ij} \in S_0$, such that $a_{ij} = 0$, $\forall i = i_s$, for $s = 1, \ldots, k$, i.e. $A$ is a matrix with $k$ fixed zero rows. It is easy to check that sum and product of two matrices of this type are also matrices of this type.

 Let $S$ be a semiring of matrices of this type. For any matrix $A = (a_{ij}) \in S$ we define

 $$\delta(A) = (a_{ij}),\; \mbox{where}\; a_{ij} = 0, \forall i = i_s\; \mbox{and}\; a_{ij} = 0,\forall j \neq i_s, \; \mbox{where}\; s = 1, \ldots, k.$$

Let $A, B \in S$, where $A = (a_{ij})$ and $B = (b_{ij})$. Since $A + B \in S$, it follows that $\delta (A + B) = (a_{ij} + b_{ij})$, where
 $a_{ij} + b_{ij} = 0, \forall i = i_s$ and  $a_{ij} + b_{ij} = 0, \forall j \neq i_s$, $s = 1, \ldots, k$. Clearly $\delta (A + B) = \delta (A) + \delta (B)$.

 Let $A\,B = C = (c_{ij}$. Then $\displaystyle c_{ij} = \sum_{l=1}^n a_{il}b_{lj}$. Since $a_{il} = 0$, $\forall i = i_s$, it follows $c_{ij} = 0$, $\forall i = i_s$, $s = 1, \ldots, k$. Now, $\delta(C) = (c_{ij})$, where  $c_{ij} = 0$, $\forall i = i_s$ and  $c_{ij} = 0$, $\forall j \neq i_s$, where $s = 1, \ldots, k$. It is easy to prove that $\delta(A)\,B = 0$. Indeed, the elements of $\delta(A)\,B$ are $\displaystyle d_{ij} = \sum_{l=1}^n a_{il}b_{lj}$. Since $a_{il} = 0$, $\forall l \neq i_s$ and $b_{lj} = 0$, $\forall l = i_s$, it follows $d_{ij} = 0$, $\forall i, j = 1, \ldots, n$.

 Let $A\, \delta(B) = (f_{ij})$, where $\displaystyle f_{ij} = \sum_{l=1}^n a_{il}b_{lj}$. Let $i = i_s$,  $s = 1, \ldots, k$. Now $a_{il} = 0$ and then $f_{ij} = 0$, where $i = i_s$. Let $j \neq i_s$,  $s = 1, \ldots, k$. Then $b_{lj} = 0$ and hence $f_{ij} = 0$, where $j \neq i_s$.

 So, $c_{ij}$ and $f_{ij}$ are sums ot the type $\sum_{l=1}^n a_{il}b_{lj}$ satisfying the same conditions. Hence $c_{ij} = f_{ij}$,  $\forall i, j = 1, \ldots, n$. Thus $\delta(A\,B) = A\, \delta(B)$ and $\delta(A\,B) = A\, \delta(B) + \delta(A)\,B$. Hence $\delta$ is a derivation.

\vspace{5mm}

\begin{center}
 \textbf{$\mathbf{(0,1)}$ -- matrices}
\end{center}

\vspace{3mm}

Let $S_0$ be an additively idempotent semiring with zero and identity element, $M_n(S_0)$ be the semiring of $n \times n$ - matrices with entries $S_0$ and $S$ is subsemiring of $M_n(S_0)$.

If $A = (a_{ij}) \in S$ the \textbf{pattern} of $A$ is called a matrix $P_A = (p_{ij})$, where
$$p_{ij} = \left\{\begin{array}{ll}1, & \; \mbox{if}\; a_{ij} \neq 0\\
0, & \; \mbox{if}\; a_{ij} = 0\\ \end{array} \right. .$$

If $B = (b_{ij}) \in S$, then $A + B = (a_{ij} + b_{ij})$. Now $P_{A+B} = (r_{ij})$, where $r_{ij} = 1$, when
$a_{ij} + b_{ij} \neq 0$ and $r_{ij} =0$ when $a_{ij} + b_{ij} = 0$. (Note that if $\alpha + \beta = 0$, where $\alpha, \beta \in S_0$, it follows $\alpha = \alpha + 0 = \alpha + \alpha + \beta = \alpha + \beta = 0$.) Hence
$$ a_{ij} + b_{ij} = 0 \iff a_{ij} = 0 \; \mbox{and}\; b_{ij} = 0.$$

If $P_B = (q_{ij})$, then $q_{ij} = 0 \iff b_{ij} = 0$. Thus $r_{ij} = p_{ij} + q_{ij} = 1$ if at least one of $a_{ij}$ and $b_{ij}$ is not zero and $r_{ij} = p_{ij} + q_{ij} = 0$ if $a_{ij} = b_{ij} = 0$. So, we prove that $P_{A + B} = P_A + P_B$.

\vspace{1mm}

Let $C = A\,B = (c_{ij})$. Then $\displaystyle c_{ij} = \sum_{k=1}^n a_{ik}b_{kj}$, $i, j = 1, \ldots, n$. Let $P_C = (s_{ij})$, where $s_{ij} = 1$ if $c_{ij} \neq 0$ and $s_{ij} = 0$ if $c_{ij} = 0$. On the other hand $P_A\,P_B = (p_{ij})(q_{ij}) = (u_{ij})$, where $\displaystyle u_{ij} = \sum_{k=1}^n p_{ik}q_{kj}$. Let us note that $u_{ij} = 1$ if at least one of products $p_{ik}q_{kj}$ is equal to 1 and $u_{ij} = 0$ if all of products $p_{ik}q_{kj}$ are equal to 0.

If  all of products $p_{ik}q_{kj}$, $k = 1, \ldots, n$, are equal to 0, then all of products $a_{ik}b_{kj}$, $k = 1, \ldots, n$, are equal to 0, so $c_{ij} = 0$ and hence $s_{ij} = 0$. If $u_{ij} = 1$, then $p_{ik}q_{kj} = 1$ for some $k$, $k = 1, \ldots, n$. Then $a_{ik}b_{kj} \neq 0$ for this $k$, which implies $c_{ij} \neq 0$ and hence $s_{ij} = 1$. So, we prove that $P_{A\,B} = P_A\, P_B$.
Thus we prove

\vspace{3mm}

\textbf{Lemma 1.} \textsl{ The set of patterns $P_A$ of matrices $A \in S$ is a subsemiring of $M_n(S_0)$.}

\vspace{2mm}

We denote this semiring by $P(S)$. If for any matrix $A \in S$ it follows $P_A \in S$, then $P(S)$ is a subsemiring of $S$. The matrix $P_B$ is called \textbf{submatrix} of the matrix $P_A$ if for some indices $i$ and $j$, where $i, j = 1, \ldots, n$, from
$b_{ij} = 1$, it follows $a_{ij} = 1$.
\vspace{2mm}

Let us consider the semiring $S = Diag_n(S_0)$ of diagonal $n \times n$ - matrices  with entries from an additively idempotent semiring $S_0$. Since the pattern $P_A$ of any matrix $A \in S$ is also an element of $S$, then $P(S)$ is a subsemiring of $S$. We have the following description of this semiring:
$$P(S) = \left\{ \sum_{i = 1}^n \varepsilon_i E_{ii}\; |\;  \varepsilon_i = 0 \; \mbox{or}\;  \varepsilon_i = 1,\; i = 1, \ldots, n\right\}.$$

\vspace{2mm}

\textsl{\textbf{Lemma 2.} Let $S$ be an additively idempotent semiring and $Z(S)$ be a center of $S$. Then for any $x \in Z(S)$ the map $\delta_x : S \rightarrow S$ such that $\delta_x(a) = xa$, where $a \in S$, is a derivation.}

\emph{Proof.} Let $a, b \in S$. Then $\delta_x(a) = xa$, $\delta_x(b) = xb$, $\delta_x(a + b) = x(a + b) = xa + xb = \delta_x(a) + \delta_x(b)$ and $\delta_x(ab) = xab  = xab + xab = \delta_x(a)b + axb = \delta_x(a)b + a\delta_x(b)$.

\vspace{3mm}

Since the additively idempotent semirings are not rings, there is not a concept commutator. So the derivations $\delta_x : S \rightarrow S$ such that $\delta_x(a) = xa$, where $a \in S$ are called \textbf{inner} derivations in additively idempotent semiring $S$.
\vspace{2mm}

Clearly for $S = Diag_n(S_0)$ the semiring $P(S)$ is a subsemiring of the center of $S$. So, for any matrix $P_X \in P(S)$ from Lemma 2 follows that the map $\delta_X : S \rightarrow S$ such that $\delta_X(A) = A\,P_X$ is an inner derivation.
\vspace{2mm}

If we define for two inner derivations $\delta_X$ and $\delta_Y$ the sum $\delta_X + \delta_Y = \delta_{X + Y}$ and the product  $\delta_X \delta_Y = \delta_{X\,Y}$ it follows that the set of these derivations is a subsemiring of $S$ isomorphic to $P(S)$. This semiring is  denoted by $D(S)$.

\vspace{5mm}

\textbf{Theorem 1.} \textsl{Let   $S = Diag_n(S_0)$ be the semiring of diagonal $n \times n$ - matrices  with entries from an additively idempotent semiring $S_0$, $P(S)$ be a subsemirings of patterns of matrices of $S$  and $D(S)$ be the semiring of the inner derivation $\delta_X$, such that $\delta_X(A) = A\,P_X$ for any matrix $P_X \in P(S)$.}

\textsl{a) The zero of  $D(S)$ is the derivation $\delta_0 =0$. The unit element of  $D(S)$ is the identity map $\delta_E = i$.}

\textsl{b) $\delta_X^2 = \delta_X$ for any derivation $\delta_X \in D(S)$.}

\textsl{c) $\delta_X \delta_Y = \delta_Y \delta_X$ for any derivations $\delta_X, \delta_Y \in D(S)$.}

\textsl{d) $\delta_{E_{ii}}\,\delta_{E_{jj}} = 0$, if $i \neq j$.}

\textsl{e) $\delta_X + \delta_Y = 0 \iff \delta_X = \delta_Y = 0$ and $\;\delta_X \,\delta_Y = i \iff \delta_X = \delta_Y = i$.}

\textsl{f) The sets} $$I_X = \{\delta_Y \; | \; P_Y \; \mbox{is a submatrix of}\; P_X\}$$
\textsl{are all the ideals of the semiring $D(S)$. Minimal ideals of $D(S)$ are $I_{E_{ii}}$, maximal ideals of $D(S)$ are $I_{E - E_{ii}}$, where $i = 1, \ldots, n$.}

\textsl{g) If $\delta$ is an arbitrary (not necessary inner) derivation in semiring $S$, then $\delta(E_{ii}) = \alpha_i E_{ii}$, where $\alpha_i \in S_0$ for all $i = 1, \ldots, n$.}
\vspace{2mm}

\emph{Proof.}  All statements from a) to f) are clear.

g) Since $E_{ii}^2 = E_{ii}$, it follows $E_{ii}\delta(E_{ii}) = \delta(E_{ii})$. But $\delta(E_{ii}) \in S$, so $\delta(E_{ii}) = \alpha_1 E_{11} + \cdots + \alpha_n E_{nn}$, where $\alpha_i \in S_0$ for $i = 1, \ldots, n$. Then $E_{ii}(\alpha_1 E_{11} + \cdots + \alpha_n E_{nn}) = \delta(E_{ii})$ and hence $\delta(E_{ii}) = \alpha_i E_{ii}$.
 \vspace{3mm}

\textbf{\emph{Remark.}} Note that $D(S)$ is a subsemiring in the semiring of inner derivations $$ID(S) = \big\{\delta^X\; |\; X = \sum_{i=1}^n c_i E_{ii}\big\},$$
 where $c_i$, $i = 1, \ldots, n$ are elements from the center of the additively idempotent semiring $S_0$ and $\delta^X(A) = AX$ for any $A \in S$. It is easy to prove that the inner derivations commute, i.e. $\delta^X \delta^Y = \delta^Y \delta^X$ for any derivations $\delta^X, \delta^Y \in ID(S)$
 \vspace{3mm}

 Are there non diagonal matrices such that their derivatives are products of matrix by some $P_X \in P(S)\,$?
  \vspace{3mm}

  \textbf{\emph{Example 3. }}  Let positive integers $i_1, \ldots, i_k$ are fixed, where $k < n$. Let us consider the matrices $A = (a_{ij})$, $i, j = 1, \ldots n$, where $a_{ij} \in S_0$, such that $a_{ij} = 0$, $\forall i = i_s$ and $\forall j = i_s$, for  $s = 1, \ldots, k$, i.e. $A$ is a matrix with $k$ fixed diagonal elements equal to zero and also all  elements of the row and the column, containing this zero of the main diagonal, are also equals to zero.  It is easy to check that sum and product of two matrices of this type are also matrices of this type.
 Let $S$ be a semiring of matrices of this type.

 Let us consider a matrix $\displaystyle P_X = \sum_{i = 1}^n \varepsilon_i E_{ii}$ such that $\varepsilon_{i_s} = 0$ for all $s = 1, \ldots, k$. Clearly $P_X \in S$. It is easy to see that $A\,P_X = P_X\,A$ for any $A \in S$ and hence $A\,P_Y = P_Y\,A$ for all submatrices $P_Y$ of $P_X$. Thus we prove that all $\delta_Y \in I_X\,$ (see definition in f) of Theorem 1) are derivation in the semiring $S$.

\vspace{5mm}

\begin{center}
 \textbf{Derivatives of triangular matrices}
\end{center}

 \vspace{3mm}

We start with two examples of derivations in upper triangular matrices.
 \vspace{3mm}

  \textbf{\emph{Example 4. }}  Let $UTM_n(S_0)$ be the semiring ot the upper triangular $n \times n$ - matrices with entries from the additively idempotent semiring $S_0$.
\vspace{1mm}

Let for $A = (a_{ij}) \in UTM_n(S_0)$ we define

 $$ d(A) = A \backslash diag(a_{11},\ldots, a_{nn}).  $$

For some $B = (b_{ij}) \in UTM_n(S_))$ the arbitrary element of the matrix  $d(AB)$ has the form $$ \sum_{k=i}^j a_{ik}b_{kj},$$ where $1 \leq i \leq k \leq j \leq n.$

On the other hand the arbitrary element of the matrix $d(A)B$ is $\displaystyle \sum_{1 < k \leq j} a_{ik}b_{kj}$.

 Similarly the arbitrary element of the matrix $Ad(B)$ is $\displaystyle \sum_{i \leq k < j} a_{ik}b_{kj}$.

 Since $S_)$ is an additively idempotent semiring, it follows that $$ \sum_{i < k < j} a_{ik}b_{kj} + \sum_{i < k < j}  a_{ik}b_{kj} = \sum_{i < k < j}  a_{ik}b_{kj}.$$ Hence $d(AB) = d(A)B + Ad(B)$ and since $d$ is evidently linear, it follows that $d$ is a derivation.
\vspace{3mm}

 Obviously the derivation $d$ is idempotent.

 \vspace{5mm}

  \textbf{\emph{Example 5. }}  Let us consider the  semiring $UTM_n(S_0)$ of the upper triangular $n \times n$ - matrices with entries in arbiytary semiring $S_0$.
For a matrix  $A = (a_{ij}) \in UTM_n(S)$ we define a two maps  by rules
$$\delta_1(A) = a_{11}E_{11} + a_{12}E_{12} + \cdots a_{1n}E_{1n}\; \mbox{and}\;    \delta_2(A) = a_{12}E_{12} + \cdots a_{1n}E_{1n}.$$

Analogously for a matrix $B = (b_{ij}) \in UTM_n(S)$ we have
$$\delta_1(B) = b_{11}E_{11} + b_{12}E_{12} + \cdots b_{1n}E_{1n}\; \mbox{and}\;    \delta_2(B) = b_{12}E_{12} + \cdots b_{1n}E_{1n}.$$

 Hence
$$\delta_1(AB) = a_{11}b_{11}E_{11} + (a_{11}b_{12} + a_{12}b_{22})E_{12} + \cdots + \sum_{j=1}^n a_{1j}b_{jn}\, E_{1n}\;$$ and $$\delta_2(AB) = (a_{11}b_{12} + a_{12}b_{22})E_{12} + \cdots + \sum_{j=1}^n a_{1j}b_{jn}\, E_{1n}.$$

Now calculate $$ \delta_1(A)B = a_{11}b_{11}E_{11} + a_{12}b_{22}E_{12} + \cdots + \sum_{j=2}^n a_{1j}b_{jn}\, E_{1n}\;\; \mbox{and} \;  A\delta_1(B) = a_{11}b_{11}E_{11} + a_{11}b_{12}E_{12} + \cdots + a_{11}b_{1n}\, E_{1n}.$$

If we now suppose that $S_0$ is an additively idempotent semiring, it follows that $\delta_1(A)B +  A\delta_1(B) = \delta_1(AB)$.

In arbitrary semiring $S_0$ we have
 $$\delta_2(A)B +  A\delta_2(B) = (a_{11}b_{12} + a_{12}b_{22})E_{12} + \cdots + \sum_{j=1}^n a_{1j}b_{jn}\, E_{1n} = \delta_2(AB).$$

Clearly $\delta_1$ and $\delta_2$ are linear maps. So, it follows that $\delta_1$ is a derivation, when $S_0$ is an additively idempotent semiring, but $\delta_2$ is a derivation, when $S_0$ is arbitrary semiring.
 \vspace{1mm}

Note that  $\delta_1$ and $\delta_2$ are idempotent derivations.

The map $\varphi : UTM_n(S_0) \rightarrow UTM_n(S_0)$, such that $\varphi(A) = a_{1k}E_{1k}  + \cdots a_{1n}E_{1n}$, where $A = (a_{ij}) \in UTM_n(S_0)$ is not a derivation when $k \geq 3$.
 \vspace{4mm}

  In the next examples we  show that there are inner derivations in subsmirings of $UTM_n(S_0)$, where $S_0$ is an additively idempotent semiring.

  \vspace{2mm}
 \textbf{\emph{Example 6. }} Let $S_n$ be a semiring of all matrices $A= (a_{ij})$, $A \in UTM_n(S_0)$, such that $a_{11} = a_{nn}$. Obviously $E_{1n} \in S_n$. Now, for $A \in S_n$ calculate
$$A\,E_{1n} = a_{11}E_{1n}\;\; \mbox{and}\;\; E_{1n}\,A = a_{nn}E_{1n}.$$

Hence $A\,E_{1n} =  E_{1n}\,A$. From  Lemma 2, it follows that $\delta_n(A) = a_{11}E_{1n}$
  define a derivation in semiring $S_n$. Note that this derivation is nilpotent since $\delta_n^2(A) = 0$.

 \vspace{2mm}

 \textbf{\emph{Example 7. }}   Now consider the matrix $A = (a_{ij}) \in UTM_n(S_0)$ such that $$a_{11} = a_{n-1n-1}, \;\; a_{12} = a_{n-1n}\;\; \mbox{and}\;\; a_{22} = a_{nn}.$$
  It is easy to show that if $A$ and $B$ are matrices of this type then also $A+B$ and $AB$ are matrices of this type. Let $S_{n-1}$ be a subsemiring of $UTM_n(S_0)$ consisting of matrices $A$ of this type. Obviously $E_{1n-1} + E_{2n} \in S_{n-1}$. For any $A \in S_{n-1}$ calculate
$$A(E_{1n-1} + E_{2n}) = a_{11}E_{1n-1} + a_{12}E_{1n} + a_{22}E_{2n} = \hphantom{aaaaaaaaaaaaaaaa}$$ $$\hphantom{aaaaaaaaaaaaa} = a_{n-1n-1}E_{1n-1} + a_{n-1n}E_{1n} + a_{nn}E_{2n} = (E_{1n-1} + E_{2n})A.$$

From   Lemma 2, it follows that $\delta(A) = a_{11}E_{1n-1} + a_{12}E_{1n} + a_{22}E_{2n}$ or
$$\delta \left(
    \begin{array}{ccccc}
      a_{11} & a_{12} & \cdots & a_{1 n-1} & a_{1n} \\
      0 & a_{22} & \cdots & a_{2 n-1} & a_{2n} \\
      \cdots & \cdots & \cdots & \cdots & \cdots \\
      0 & 0 & \cdots & a_{11} & a_{12} \\
       0 & 0 & \cdots & 0 & a_{22} \\
    \end{array}
  \right) = \left(
    \begin{array}{ccccc}
      0 & 0 & \cdots & a_{11} & a_{12} \\
      0 & 0 & \cdots & 0  & a_{22}\\
      \cdots & \cdots & \cdots & \cdots & \cdots \\
      0 & 0 & \cdots & 0 & 0 \\
      0 & 0 & \cdots & 0 & 0 \\
    \end{array}
  \right)$$
  define a derivation in semiring $S_{n-1}$. Note that this derivation is also nilpotent.
\vspace{3mm}

In the last step of search for an inner derivation in subsemiring of $UTM_n(S_0)$ we consider the matrix
$$D = E_{12} + E_{23} + \cdots + E_{n-1\,n}.$$

Let $A = (a_{ij}) \in UTM_n(S_0)$, that is $\displaystyle A = \sum_{i \leq j} a_{ij}E_{ij}$, where $i,j = 1, \ldots n$. We shall find the conditions of elements of the matrix $A$ such that $A\,D = D\,A$. Using the equalities
$$E_{ik}\,E_{lj} = \left\{ \begin{array}{rl} E_{ij}, & \mbox{if}\;\; k = l,\\ 0, & \mbox{if}\;\; k \neq l \end{array} \right.$$
we calculate
$$A\,D = \left(\sum_{i \leq j} a_{ij}E_{ij}\right)(E_{12} + E_{23} + \cdots + E_{n-1\,n}) = \sum_{i \leq j} a_{ij}E_{ij+1},$$
where $i = 1, \ldots n-1$ and $j = 2, \ldots, n$ or
$$A\,D = \left(
    \begin{array}{cccccc}
      0 & a_{11} & a_{12} & \cdots & a_{1\,n-2} & a_{1\,n-1} \\ \vspace{1mm}
      0 & 0 & a_{22} &\cdots & a_{2\,n-2}  & a_{2\,n-1}\\ \vspace{1mm}
      \vdots & \vdots& \vdots & \vdots & \vdots & \vdots \\
      0 & 0 & 0& \cdots & 0 & a_{n-1\,n-1} \\
      0 & 0 & 0 & \cdots & 0 & 0 \\
    \end{array}
  \right).$$

  Using  same reasoning we calculate $$D\,A = (E_{12} + E_{23} + \cdots + E_{n-1\,n})\left(\sum_{i \leq j} a_{ij}E_{ij}\right) = \sum_{i \leq j} a_{ij}E_{i-1j},$$
where $i = 2, \ldots n$ and $j = 2, \ldots, n$ or
$$D\,A = \left(
    \begin{array}{cccccc}
      0 & a_{22} & a_{23} & \cdots & a_{2\,n-1} & a_{2\,n} \\ \vspace{1mm}
      0 & 0 & a_{33} &\cdots & a_{3\,n-1}  & a_{3\,n}\\ \vspace{1mm}
      \vdots & \vdots& \vdots & \vdots & \vdots & \vdots \\
      0 & 0 & 0& \cdots & 0 & a_{n\,n} \\
      0 & 0 & 0 & \cdots & 0 & 0 \\
    \end{array}
  \right).$$

  Then by comparing the elements of $AD$ and $DA$ we find:
  $$\begin{array}{c}a_{11} = a_{22} = \cdots = a_{n-1\,n-1} = a_{nn},\\
  a_{12} = a_{23} = \cdots = a_{n-1\,n}\\
  \cdots \cdots \cdots \cdots \cdots \cdots \cdots\\
  a_{1\,n-2} = a_{2\,n-1} = a_{3\,n},\\
  a_{1\,n-1} = a_{2\,n}\end{array} .$$

  Thus we prove

\vspace{3mm}

\textbf{Proposition 2.} \textsl{ Let $A \in UTM_n(S_0)$, where $S_0$ is an additively idempotent semiring and $D = E_{12} + E_{23} + \cdots + E_{n-1\,n}$. Then $A\,D = D\,A$ if and only if the matrix $A$ has the form}
$$\left(
    \begin{array}{cccccc} \vspace{1.5mm}
      a_{0} & a_{1} & a_2 &\cdots & a_{n-2} & a_{n-1} \\ \vspace{1.5mm}
      0 & a_{0} & a_1 & \cdots & a_{n-3} & a_{n-2} \\
      0 & 0 & a_{0} & \cdots & a_{n-4} & a_{n-3}\\
      \vdots & \vdots & \vdots & \vdots & \vdots & \vdots \\ \vspace{1.5mm}
      0 & 0 & 0& \cdots & a_{0} & a_{1} \\
       0 & 0 & 0 &\cdots & 0 & a_{0} \\
    \end{array}\right). \eqno{(1)}$$

     \vspace{3mm}

\textbf{\emph{Remark.}} Clearly that if $A$ has the form (1) then $A\,D = D\,A$.

\vspace{5mm}

\begin{center}
 \textbf{Derivatives of Toeplitz matrices}
\end{center}

 \vspace{3mm}

In linear algebra is used \textbf{diagonal-constant matrices} or \textbf{Toeplitz matrices} named after \emph{Otto Toeplitz} (1881 - 1940) -- german mathematician working in functional analysis.
A Toeplitz matrix is an $n \times n$ matrix $A = (a_{ij})$, where $i,j = 0, \ldots, n-1$, $a_{ij} \in \mathbb{R}$  and $a_{ij} = a_{j-i}$, i.e. a matrix of the form
$$A = \left(
    \begin{array}{cccccc} \vspace{1.5mm}
      a_{0} & a_{1} & a_2 &\cdots & a_{n-2} & a_{n-1} \\ \vspace{1.5mm}
      a_{-1} & a_{0} & a_1 & \cdots & a_{n-3} & a_{n-2} \\
      a_{-2} & a_{-1} & a_{0} & \cdots & a_{n-4} & a_{n-3}\\
      \vdots & \vdots & \vdots & \vdots & \vdots & \vdots \\
      a_{-(n-2)} & a_{-(n-3)} & a_{-(n-4)}& \cdots & a_{0} & a_{1} \\
       a_{-(n-1)} & a_{-(n-2)} & a_{-(n-3)} &\cdots & a_{-1} & a_{0} \\
    \end{array}\right).$$
 Such matrices arise in applications including areas of coding, spectral estimation, harmonic analysis, adaptive filtering, graphical models, noise reduction. Toeplitz matrices also appears in problems  in physics, statistics and signal processing -- see [7].

\vspace{3mm}

\textbf{Proposition 3.} \textsl{ Let $S_0$ be an additively idempotent semiring, $A$ is a Toeplitz matrix with entries from $S_0$ and $D = E_{12} + E_{23} + \cdots + E_{n-1\,n}$. Then $A\,D = D\,A$ if and only if the matrix $A$ has the form (1)}
\vspace{2mm}

\emph{Proof.} The proof is similar to the proof of Proposition 2.
\vspace{3mm}

     Let us consider the upper triangular Toeplitz matrices, i.e. matrices of the form (1) with entries from additively idempotent semiring $S_0$.

It is easy to see that all such matrices form an additively idempotent semiring. This semiring  we denote by $\mathbb{UT}_n(S_0)$.

For the matrix $D = E_{12} + \cdots + E_{i-1 i} + \cdots + E_{n-1 n} \in \mathbb{UT}_n(S_0)$, using  $E_{i-1 i} \cdot E_{i i+1} = E_{i-1 i + 1}$, it follows
$$\begin{array}{l} D^2 = E_{13} + E_{24} + \cdots + E_{n-2 n}\\
D^3 = E_{14} + E_{25} + \cdots + E_{n-3 n}\\
\cdots \;\; \cdots\\
D^{n-1} = E_{1 n}\\
D^n = 0
\end{array}.$$

Hence any matrix $A \in \mathbb{UT}_n(S_0)$ has a form
$$A = a_0 E + a_1 D + \cdots + a_{n-1} D^{n-1},$$
where $a_i \in S_0$, $i = 0, 1, \ldots, n-1$.

So, the semiring of patterns of matrices in $\mathbb{UT}_n(S_0)$ is
$$P\left(\mathbb{UT}_n(S_0)\right) = \left\{\varepsilon_0 E + \varepsilon_1 D + \cdots + \varepsilon_{n-1} D^{n-1} \; |\;  \varepsilon_i = 0 \; \mbox{or}\;  \varepsilon_i = 1,\; i = 0, \ldots, n - 1\right\}.$$

If $X \in P\left(\mathbb{UT}_n(S_0)\right)$ from  Lemma 2, it follows that $\delta_X(A) = AX$ is a derivation  in semiring $\mathbb{UT}_n(S_0)$. Like in the case of diagonal matrices we show that these inner derivations form a semiring isomorphic to $P\left(\mathbb{UT}_n(S_0)\right)$. This semiring is  denoted by $D\left(\mathbb{UT}_n(S_0)\right)$.

The following result is similar to Theorem 1.
\vspace{4mm}

\textbf{Theorem 2.} \textsl{Let   $\mathbb{UT}_n(S_0)$  be the semiring of upper Toeplitz  $n \times n$ - matrices  with entries from an additively idempotent semiring $S_0$, $P\left(\mathbb{UT}_n(S_0)\right)$ be a subsemiring of patterns of matrices of $\mathbb{UT}_n(S_0)$  and $D\left(\mathbb{UT}_n(S_0)\right)$ be the semiring of the inner derivations $\delta_X$, such that $\delta_X(A) = A\,P_X$ for any matrix $P_X \in P\left(\mathbb{UT}_n(S_0)\right)$.}

\textsl{a) The zero of  $D\left(\mathbb{UT}_n(S_0)\right)$ is the derivation $\delta_0 =0$. The unit element of  $D\left(\mathbb{UT}_n(S_0)\right)$ is the identity map $\delta_E = i$.}

\textsl{b) $\delta_X^2 = \delta_X$ for any derivation $\delta_X \in D\left(\mathbb{UT}_n(S_0)\right)$ such that the pattern matrix $P_X$ has the form $E +  D^k + \cdots + \varepsilon_{n-1} D^{n-1}$, where $2k > n-1$.}

\textsl{c) $\delta_{X}$ is a nilpotent derivation if and only if the pattern matrix $P_X$ has the form $\varepsilon_1 D + \cdots + \varepsilon_{n-1} D^{n-1}$.}

\textsl{d) $\delta_X \delta_Y = \delta_Y \delta_X$ for any derivations $\delta_X, \delta_Y \in D\left(\mathbb{UT}_n(S_0)\right)$.}

\textsl{e) $\delta_X + \delta_Y = 0 \iff \delta_X = \delta_Y = 0$ and $\;\delta_X \,\delta_Y = i \iff \delta_X = \delta_Y = i$.}

\textsl{f) The set} $$I_D = \{  \alpha_1 D + \cdots + \alpha_{n-1} D^{n-1}\; | \; \alpha_i = S_0,\; i = 1, \ldots, n - 1\}$$
\textsl{is the maximal ideal of the semiring $\mathbb{UT}_n(S_0)$. }

\textsl{g) The ideal $I_D$ is closed under the arbitrary derivation $\delta$ in the semiring $\mathbb{UT}_n(S_0)$.}
\vspace{2mm}

\emph{Proof.}  It is easy to see that all statements from a) to f) are satisfied.

g) Since $D^n = 0$, it follows $\delta(D)\,D^{n-1} = 0$. But $\delta(D) \in \mathbb{UT}_n(S_0)$, hence $\delta(D) = a_0 E + a_1 D + \cdots + a_{n-1} D^{n-1}$, where $a_i = S_0$  for  $i = 0, \ldots, n - 1$. Then $0 = \delta(D)\,D^{n-1} = a_0\,D^{n-1}$ which implies $a_0 = 0$. Thus we prove that $\delta(D) \in I_D$. Since $\delta(D^k) = \delta(D)\,D^{k-1}$ for any $k = 2, \ldots, n-1$, it follows that the ideal $I_D$ is closed under the arbitrary derivation $\delta : \mathbb{UT}_n(S_0) \rightarrow \mathbb{UT}_n(S_0)$.

\vspace{3mm}

\textbf{\emph{Remark.}} Note that $D\left(\mathbb{UT}_n(S_0)\right)$ is a subsemiring in the semiring of inner derivations $$ID\left(\mathbb{UT}_n(S_0)\right) = \big\{\delta^X\; |\; X = c_0 E + c_1 D + \cdots + c_{n-1} D^{n-1}\big\},$$
 where $c_i$, $i = 0, \ldots, n-1$ are elements from the center of the additively idempotent semiring $S_0$ and $\delta^X(A) = AX$ for any $A \in \mathbb{UT}_n(S_0)$. It is easy to prove that the arbitrary inner derivation $\delta^X$ is nilpotent when $X = c_1 D + \cdots + c_{n-1} D^{n-1}$, where $c_i$ are from the center of $S_0$ and that  the inner derivations commute, i.e. $\delta^X \delta^Y = \delta^Y \delta^X$ for any derivations $\delta^X, \delta^Y \in ID\left(\mathbb{UT}_n(S_0)\right)$.
\vspace{3mm}

Now consider three maps in semiring $\mathbb{UT}_4(S_0)$.

First map
$$d\,\big(\;\left(
     \begin{array}{cccc}
       a_0 & a_1 & a_2 & a_3 \\
       0 & a_0 & a_1 & a_2 \\
       0 & 0 & a_0 & a_1 \\
       0 & 0 & 0 & a_0 \\
     \end{array}
   \right)\;\big) = \left(
     \begin{array}{cccc}
       0 & a_1 & a_2 & a_3 \\
       0 & 0 & a_1 & a_2 \\
       0 & 0 & 0 & a_1 \\
       0 & 0 & 0 & 0 \\
     \end{array}
   \right)$$

is a derivation, defined for arbitrary triangular matrices -- see Example 4.

Second map
$$\delta\,\big(\;\left(
     \begin{array}{cccc}
       a_0 & a_1 & a_2 & a_3 \\
       0 & a_0 & a_1 & a_2 \\
       0 & 0 & a_0 & a_1 \\
       0 & 0 & 0 & a_0 \\
     \end{array}
   \right)\;\big) = \left(
     \begin{array}{cccc}
       0 & 0 & 0 & a_3 \\
       0 & 0 & 0 & 0 \\
       0 & 0 & 0 & 0 \\
       0 & 0 & 0 & 0 \\
     \end{array}
   \right)$$

is a derivation, defined for  triangular matrices of special type -- see Example 6.

Third  map
$$\varphi\,\big(\;\left(
     \begin{array}{cccc}
       a_0 & a_1 & a_2 & a_3 \\
       0 & a_0 & a_1 & a_2 \\
       0 & 0 & a_0 & a_1 \\
       0 & 0 & 0 & a_0 \\
     \end{array}
   \right)\;\big) = \left(
     \begin{array}{cccc}
       0 & 0 & a_2 & a_3 \\
       0 & 0 & 0 & a_2 \\
       0 & 0 & 0 & 0 \\
       0 & 0 & 0 & 0 \\
     \end{array}
   \right)$$
is \emph{{not}} a derivation.

In order to find some conditions of the matrices such that the map $\varphi$ to become a derivation we prove the following proposition:
\vspace{3mm}

\textsl{\textbf{ Proposition 4.} Let matrix $A \in \mathbb{UT}_n(S_0)$, where $S_0$ is an arbitrary semiring,  has a form
 $$A = a_0 E + a_{n-k} D^{n-k} + \cdots + a_{n-1} D^{n-1},$$
 where $1 \leq k \leq n-1$.  The set of these matrices is a semiring $S_k^n$ and the map $\delta : S_k^n \rightarrow S_k^n$ such that
 $$\delta(A) =  a_{n-k} D^{n-k} + \cdots + a_{n-1} D^{n-1}$$ is a derivation. For $k + 1 \leq n < 2k$ the semiring $S_0$ must be an additively idempotent, but for $2k \leq n$ this condition may be ommited.}
\vspace{2mm}

\emph{Proof.} For $B = b_0 E + b_{n-k} D^{n-k} + \cdots + b_{n-1} D^{n-1} \in \mathbb{UT}_n(S_0)$ we have $\delta(B) = b_{n-k} D^{n-k} + \cdots + b_{n-1} D^{n-1}$. It is clear that $A + B$ and $AB$ are matrices of the same form, so the set of these matrices $S_k^n$ is a subsemiring of $\mathbb{UT}_n(S_0)$.
\vspace{1mm}

\emph{Case 1.} Let $2k \leq n$. Since
$$AB = a_0b_0 E + (a_0b_{n-k} + b_0a_{n-k}) D^{n-k} + \cdots + (a_0b_{n-1} + b_0a_{n-1}) D^{n-1},$$
it follows $\delta(AB) = (a_0b_{n-k} + b_0a_{n-k}) D^{n-k} + \cdots + (a_0b_{n-1} + b_0a_{n-1}) D^{n-1}$. But
$$A\delta(B) = a_0b_{n-k}D^{n-k} + \cdots + a_0b_{n-1}D^{n-1}\;\;\mbox{and}\;\; \delta(A)B = b_0a_{n-k}D^{n-k} + \cdots + b_0a_{n-1}D^{n-1}.$$
Hence $\delta(AB) = A\delta(B) + \delta(A)B$. Evidently $\delta$ is a linear map, so $\delta : S_k^n \rightarrow S_k^n$ is a derivation.
\vspace{1mm}

\emph{Case 2.} Let $k + 1 \leq n < 2k$ and $2(n-k) = m$, where $m \leq n-1$. Now
$$AB = a_0b_0 E + (a_0b_{n-k} + b_0a_{n-k}) D^{n-k} + \cdots + (a_0b_{m-1} + b_0a_{m-1}) D^{m-1}  +
 (a_0b_{m} + b_0a_{m} + a_{n-k}b_{n-k}) D^{m} + $$ $$ + (a_0b_{m+1} + b_0a_{m+1} + a_{n-k}b_{n-k+1}+ a_{n-k+1}b_{n-k}) D^{m+1} + \cdots + (a_0b_{n-1} + b_0a_{n-1} + a_{n-k}b_{k-1} + \cdots + b_{n-k}a_{k-1}) D^{n-1}$$ and then
$$\delta(AB) = (a_0b_{n-k} + b_0a_{n-k}) D^{n-k} + \cdots + (a_0b_{m-1} + b_0a_{m-1}) D^{m-1}  +
 (a_0b_{m} + b_0a_{m} + a_{n-k}b_{n-k}) D^{m} + $$ $$ + (a_0b_{m+1} + b_0a_{m+1} + a_{n-k}b_{n-k+1}+ a_{n-k+1}b_{n-k}) D^{m+1} + \cdots + (a_0b_{n-1} + b_0a_{n-1} + a_{n-k}b_{k-1} + \cdots + b_{n-k}a_{k-1}) D^{n-1}.$$

We calculate
$$A\delta(B) = a_0b_{n-k}D^{n-k} + \cdots + a_0b_{m-1}D^{m-1}  +
 (a_0b_{m} +  a_{n-k}b_{n-k}) D^{m} + + (a_0b_{m+1}  + a_{n-k}b_{n-k+1}+ a_{n-k+1}b_{n-k}) D^{m+1} +$$ $$+ \cdots + (a_0b_{n-1}  + a_{n-k}b_{k-1} + \cdots + b_{n-k}a_{k-1}) D^{n-1}$$
and
$$\delta(A)B = b_0a_{n-k}D^{n-k} + \cdots + b_0a_{m-1}D^{m-1}  +
 (b_0a_{m} +  b_{n-k}a_{n-k}) D^{m} + + (b_0a_{m+1}  + b_{n-k}a_{n-k+1}+ b_{n-k+1}a_{n-k}) D^{m+1} +$$ $$+ \cdots + (b_0a_{n-1}  + b_{n-k}a_{k-1} + \cdots + a_{n-k}b_{k-1}) D^{n-1}.$$

Since $S_0$ is an additively idempotent semiring, it follows $\delta(AB) = A\delta(B) + \delta(A)B$. It is clear that $\delta$ is a linear map and then $\delta : S_k^n \rightarrow S_k^n$ is a derivation.

\vspace{5mm}

\begin{center}
 \textbf{Derivatives of circulant matrices}
\end{center}

 \vspace{3mm}

A circulant matrix is a special form of Toeplitz matrix where each row  is rotated one element to the right relative to the preceding row, i.e. a matrix of the form
$$A = \left(
    \begin{array}{ccccc} \vspace{1.5mm}
      a_{1} & a_{2} & a_3 &\cdots  & a_{n} \\ \vspace{1.5mm}
      a_{n} & a_{1} & a_2 & \cdots  & a_{n-1} \\
      a_{n-1} & a_{n} & a_{1} & \cdots  & a_{n-2}\\
      \vdots & \vdots & \vdots & \vdots  & \vdots \\
      a_{2} & a_{3} & a_{4}& \cdots  & a_{1} \\
          \end{array}\right). \eqno{(2)}$$

Circulant matrices have many applications in physics,  image processing,  probability and statistics,  numerical analysis,  number theory and cryptography -- [1] and [7].

\vspace{3mm}

\textbf{Proposition 5.} \textsl{ Let $S_0$ be an arbitrary semiring, $A$ is an arbitrary $n\times n$ -  matrix with entries from $S_0$ and $d = E_{12} + \cdots + E_{n-1\,n} + E_{n\,1}$. Then $A\,d = d\,A$ if and only if the matrix $A$ has the form (2)}
\vspace{1mm}

\emph{Proof.} For $A = (a_{ij}), i, j = 1, \ldots, n$ we calculate
$$A\,d = \left(\sum_{i =1}^n\sum_{j=1}^n a_{ij}E_{ij}\right)(E_{12} +  \cdots + E_{n-1\,n} + E_{n\,1}) = \sum_{i=1}^n\sum_{j=1}^{n-1}  a_{ij}E_{ij+1} + \sum_{k=1}^n a_{kn} E_{k1}$$
 or
$$A\,d = \left(
    \begin{array}{ccccc}
      a_{1n} & a_{11} & a_{12} & \cdots  & a_{1\,n-1} \\ \vspace{1mm}
      a_{2n} & a_{21} & a_{22} &\cdots  & a_{2\,n-1}\\ \vspace{1mm}
      \vdots & \vdots& \vdots & \vdots & \vdots \\
       a_{n-1\,n} & a_{n-1\,1} & a_{n-1\,2}& \cdots  & a_{n-1\,n-1} \\
      a_{nn} & a_{n1} & a_{n2}& \cdots  & a_{n\,n-1} \\
          \end{array}
  \right).$$

  Using  same reasoning we calculate $$d\,A = (E_{12} +  \cdots + E_{n-1\,n} + E_{n\,1})\left(\sum_{i =1}^n\sum_{j=1}^n a_{ij}E_{ij}\right) = \sum_{i=1}^{n-1}\sum_{j=1}^{n}  a_{i+1j}E_{ij} + \sum_{k=1}^n a_{1k} E_{nk}$$
 or
$$d\,A = \left(
    \begin{array}{ccccc}
      a_{21} & a_{22} & a_{23} & \cdots  & a_{2\,n} \\ \vspace{1mm}
      a_{31} & a_{32} & a_{33} &\cdots  & a_{3\,n}\\ \vspace{1mm}
      \vdots & \vdots& \vdots & \vdots  & \vdots \\
      a_{n1} & a_{n2} & a_{n3} & \cdots  & a_{n\,n} \\
      a_{11} & a_{12} & a_{13} & \cdots  & a_{1n} \\
    \end{array}
  \right).$$

  Now we  compare the elements of $A\,d$ and $d\,A$ and find:
  $$\begin{array}{ccccc} a_{1n} = a_{21}\; & a_{2n} = a_{31}\; & \cdots \; & a_{n-1\,n} = a_{n1}\; & a_{nn} = a_{11}\\
  a_{11} = a_{22}\; & a_{21} = a_{32}\; & \cdots \; & a_{n-1\,1} = a_{n2}\; & a_{n1} = a_{12}\\
  a_{12} = a_{23}\; & a_{22} = a_{33}\; & \cdots \; & a_{n-1\,2} = a_{n3}\; & a_{n2} = a_{13}\\
  \vdots &  \vdots &  \vdots &  \vdots &  \vdots \\
  a_{1\,n-1} = a_{2n}\; & a_{2\,n-1} = a_{3n}\; & \cdots \; & a_{n-1\,n-1} = a_{nn}\; & a_{n\,n-1} = a_{1n}\\
  \end{array}. $$

  If we replace $a_{1i} = a_i$ for $i = 1, \ldots, n$ we obtain that the matrix $A$ has the form (2).

  Since all circulant matrices commutes, it follows that if $A$ is a circulant matrix then $A\,d = d\,A$.
  \vspace{2mm}

  By  similar reasoning as in the last proof and using notations from Proposition 5, we prove the following

\vspace{3mm}

\textbf{Proposition 6.} \textsl{The circulant matrix $A$ commutes with the matrix $D = E_{12} + \cdots + E_{n-1\,n}$ if and only if the matrix $A = \alpha E$, where $\alpha \in S_0$}.
 \vspace{2mm}

 Since sum and product of two circulant matrices are also circulant matrices, it follows that the set of $n\times n$ - circulant matrices with entries from an additively idempotent semiring $S_0$ form a semiring which we denote by $\mathbb{CM}_n(S_0)$, The elements of this semiring are matrix $d$ and
 $$d^2 = E_{13} + E_{24} + \cdots + E_{n-2\,n} + E_{n-1\,1} + E_{n\,2},$$
  $$d^3 = E_{14} + E_{25} + \cdots + E_{n-3\,n} + E_{n-2\,1} + E_{n-1\,2} + E_{n-3\,n},$$
  $$ \cdots  \cdots   \cdots   \cdots  \cdots \cdots  \cdots \cdots \cdots \cdots $$
  $$d^{n-1} = E_{1n} + E_{21} + \cdots + E_{n\,n-1},$$
  $$d^n = E.$$

Hence any matrix $A \in \mathbb{CM}_n(S_0)$ ( i.e.of the form (2)) can be written
$$A = a_1 E + a_2 d + \cdots + a_{n} d^{n-1},$$
where $a_i \in S_0$, $i = 1, \ldots, n$.

So, the semiring of patterns of matrices in $\mathbb{CM}_n(S_0)$ is
$$P\left(\mathbb{CM}_n(S_0)\right) = \left\{\varepsilon_1 E + \varepsilon_2 d + \cdots + \varepsilon_{n} d^{n-1} \; |\;  \varepsilon_i = 0 \; \mbox{or}\;  \varepsilon_i = 1,\; i = 1, \ldots, n\right\}.$$

If $P_X \in P\left(\mathbb{CM}_n(S_0)\right)$ from  Lemma 2, it follows that $\delta_X(A) = AP_X$ is a derivation  in semiring $\mathbb{CM}_n(S_0)$. Like in the cases of diagonal and Toeplitz matrices we show that these inner derivations form a semiring isomorphic to $P\left(\mathbb{CM}_n(S_0)\right)$. This semiring is  denoted by $D\left(\mathbb{CM}_n(S_0)\right)$.

The following result is similar to Theorem 1.
\vspace{4mm}

\textbf{Theorem 3.} \textsl{Let   $\mathbb{CM}_n(S_0)$  be the semiring of circulant  $n \times n$ - matrices  with entries from an additively idempotent semiring $S_0$, $P\left(\mathbb{CM}_n(S_0)\right)$ be a subsemiring of patterns of matrices of $\mathbb{CM}_n(S_0)$  and $D\left(\mathbb{CM}_n(S_0)\right)$ be the semiring of the inner derivations $\delta_X$, such that $\delta_X(A) = A\,P_X$ for any matrix $P_X \in P\left(\mathbb{CM}_n(S_0))\right)$.}

\textsl{a) The zero of  $D\left(\mathbb{CM}_n(S_0)\right)$ is the derivation $\delta_0 =0$. The unit element of  $D\left(\mathbb{CM}_n(S_0)\right)$ is the identity map $\delta_E = i$.}

\textsl{b) Any derivation $\delta \in D\left(\mathbb{CM}_n(S_0)\right)$ is neither idempotent, nor nilpotent.}

\textsl{c) $\delta_X \delta_Y = \delta_Y \delta_X$ for any derivations $\delta_X, \delta_Y \in  D\left(\mathbb{CM}_n(S_0)\right)$.}

\textsl{d) $\delta_X + \delta_Y = 0 \iff \delta_X = \delta_Y = 0$ and $\;\delta_X  + \delta_Y = i \iff \delta_X = \delta_Y = i$.}

\textsl{e) The invertible elements of semiring $P\left(\mathbb{CM}_n(S_0)\right)$ form a cyclic group $U_P = \{E,d, \ldots, d^{n-1}\}$. The invertible elements of semiring $D\left(\mathbb{CM}_n(S_0)\right)$ form a cyclic group $U_D = \{i,\delta_d, \ldots, \delta_{d^{n-1}}\}$.}

\textsl{f) The sets} $$I_P = P\left(\mathbb{CM}_n(S_0)\right) \backslash U_P \; \; \mbox{and}\;\;    I_D = D\left(\mathbb{CM}_n(S_0)\right) \backslash U_D $$
\textsl{are maximal ideals of the semirings $P\left(\mathbb{CM}_n(S_0)\right)$ and $D\left(\mathbb{CM}_n(S_0)\right)$, respectively. }

\textsl{g) The ideal $I_P$ is closed under the arbitrary derivation $\delta$ in the semiring $P\left(\mathbb{CM}_n(S_0)\right)$.}
\vspace{2mm}

\emph{Proof.}  It is easy to see that all statements from a) to f) are satisfied.

g) Since $d^n = E$, it follows $\delta(d^n) = d^{n-1}\delta(d) = \delta(E)$. Hence $\delta(d) = d\delta(E)$ and  $\delta(d^k) = d^k\delta(E)$, where $k = 2, \ldots, n-1$. So, for any $X \in I_P$, it follows $\delta(X) = X\delta(E)$ and since $\delta(E) \in P\left(\mathbb{CM}_n(S_0)\right)$ we find that $\delta(X) \in I_P$.

\vspace{3mm}

\textbf{\emph{Remark.}} Note that $D\left(\mathbb{CM}_n(S_0)\right)$ is a subsemiring in the semiring of inner derivations $$ID\left(\mathbb{CM}_n(S_0)\right) = \big\{\delta^X\; |\; X = c_0 E + c_1 D + \cdots + c_{n-1} D^{n-1}\big\},$$
 where $c_i$, $i = 0, \ldots, n-1$ are elements from the center of the additively idempotent semiring $S_0$ and $\delta^X(A) = AX$ for any $A \in \mathbb{CM}_n(S_0)$. It is easy to prove that  the inner derivations commute, i.e. $\delta^X \delta^Y = \delta^Y \delta^X$ for any derivations $\delta^X, \delta^Y \in ID\left(\mathbb{CM}_n(S_0)\right)$.

\vspace{3mm}

Similar reasoning can be applied to block matrices and this will be considered in  next work.

\begin{center} \textbf{REFERENCES}
\end{center}

[1]
 R. Aldrovandi, Special matrices of mathematical physics (Stochastic, Circulant and Bell matrices), World Scient. Publ., 2001.

[2]
R. B. Bapat, T. E. S. Raghavan, Nonnegative matrices and applications, Cambridge University Press, 1997.

[3] R. Brualdi, H. Ryser, Combinatorial matrix theory, Cambridge Univ. Press, 1991.

[4]  P. Butkovi$\check{c}$, Max-linear Systems. Springer Monographs in
Mathematics. Springer, London, 2010.

[5]  S. I. Dimitrov, Derivations on Semirings, (to appear in this Conference proceeding).

[6]  J. S. Golan,  Semirings and Affine Equations over Them: Theory and Applications,
Springer -- Kluwer Acad.Publ., Dodrecht, (2003).

[7] R. Gray, Toeplitz and Circulant Matrices: A review., Stanford Univ. Press, 2006.

[8] T. Kollo, D. van Rosen, Advanced Multivariate Statistics with Matrices, Springer, 2005.

  [9] N. Krivulin, Complete solutions of an optimization problems in tropical semifield, arXiv:1706.00643v1 [math.OC] 2 Jun 2017.

  [10] D. Vladeva, Projections of $k$ - simplex onto the subsimplicess of arbitrary type are derivations, arXiv:1706.00033v1 [math.RA] 31 May 2017.

\end{document}